\documentclass{amsart}

\usepackage{hyperref}
\usepackage{enumerate}
\usepackage{mathrsfs}
\usepackage[matrix,arrow,graph,curve,frame]{xy}

\theoremstyle{plain}      
    \newtheorem{theorem}{Theorem}[section]
    \newtheorem{proposition}[theorem]{Proposition}

\theoremstyle{definition}

\theoremstyle{remark}
    \newtheorem*{remark}{Remark}

\newcommand{\A}{\ensuremath{\mathscr{A}}}
\newcommand{\B}{\ensuremath{\mathscr{B}}}
\newcommand{\D}{\ensuremath{\mathscr{D}}}
\newcommand{\V}{\ensuremath{\mathscr{V}}}
\newcommand{\X}{\ensuremath{\mathscr{X}}}
\newcommand{\Y}{\ensuremath{\mathscr{Y}}}
\newcommand{\Z}{\ensuremath{\mathcal{Z}}}
\newcommand{\Tamb}{\ensuremath{\mathrm{Tamb}}}

\newcommand{\LTamb}{\ensuremath{\mathrm{LTamb}}}
\newcommand{\RTamb}{\ensuremath{\mathrm{RTamb}}}
\newcommand{\VMod}{\ensuremath{\V$-$\mathbf{Mod}}}
\newcommand{\Set}{\ensuremath{\mathbf{Set}}}

\newcommand{\ox}{\ensuremath{\otimes}}
\newcommand{\op}{\ensuremath{\mathrm{op}}}
\newcommand{\pr}{\ensuremath{\mathrm{pr}}}
\newcommand{\copr}{\ensuremath{\mathrm{copr}}}

\newcommand{\ra}{\ensuremath{\xymatrix@1@=20pt{\ar[r]&}}}
\newcommand{\la}{\ensuremath{\xymatrix@1@=20pt{& \ar[l]}}}
\newcommand{\mra}{\ensuremath{\xymatrix@1@=20pt{\ar[r] |<(0.4){\object@{|}}&}}}

\begin{document}

\title{Doubles for monoidal categories}
\author{Craig Pastro and Ross Street}
\email{\{craig, street\}@maths.mq.edu.au}
\address{Centre of Australian Category Theory \\
         Department of Mathematics \\
         Macquarie University \\
         NSW 2109 Australia}
\date{\today}
\keywords{monoidal centre, Drinfeld double, monoidal category, Day
convolution}
\thanks{The first author gratefully acknowledges support of an international
Macquarie University Research Scholarship while the second gratefully
acknowledges support of the Australian Research Council Discovery Grant
DP0771252.}
\dedicatory{Dedicated to Walter Tholen on his 60th birthday}

\begin{abstract}
In a recent paper, Daisuke Tambara defined two-sided actions on an endomodule
(= endodistributor) of a monoidal $\V$-category $\A$. When $\A$ is autonomous
(= rigid = compact), he showed that the $\V$-category (that we call
$\Tamb(\A)$) of so-equipped endomodules (that we call Tambara modules) is
equivalent to the monoidal centre $\Z[\A,\V]$ of the convolution monoidal
$\V$-category $[\A,\V]$. Our paper extends these ideas somewhat. For general
$\A$, we construct a promonoidal $\V$-category $\D\A$ (which we suggest
should be called the double of $\A$) with an equivalence $[\D\A,\V] \simeq
\Tamb(\A)$. When $\A$ is closed, we define strong (respectively, left strong)
Tambara modules and show that these constitute a $\V$-category $\Tamb_s(\A)$
(respectively, $\Tamb_{l s}(\A)$) which is equivalent to the centre
(respectively, lax centre) of $[\A,\V]$. We construct localizations $\D_s \A$
and $\D_{l s} \A$ of $\D\A$ such that there are equivalences $\Tamb_s(\A)
\simeq [\D_s \A, \V]$ and $\Tamb_{l s}(\A) \simeq [\D_{l s} \A, \V]$. When
$\A$ is autonomous, every Tambara module is strong; this implies an
equivalence $\Z[\A, \V] \simeq [\D\A,\V]$.
\end{abstract}

\maketitle

\section{Introduction}

For $\V$-categories $\A$ and $\B$, a \emph{module} $T:\A \mra \B$ (also
called ``bimodule'', ``profunctor'', and ``distributor'') is a $\V$-functor
$T:\B^\op \ox \A \ra \V$. For a monoidal $\V$-category $\A$,
Tambara~\cite{Tamb} defined two-sided actions $\alpha$ of $\A$ on an
endomodule $T:\A \mra \A$. When $\A$ is autonomous (also called
``rigid'' or ``compact'') he showed that the $\V$-category $\Tamb(\A)$ of
Tambara modules $(T,\alpha)$ is equivalent to the monoidal centre
$\Z[\A,\V]$ of the convolution monoidal $\V$-category $[\A,\V]$.

Our paper extends these ideas in four ways:
\begin{enumerate}
\item our base monoidal category $\V$ is quite general (as
in~\cite{KellyBook}) not just vector spaces;

\item our results are mainly for a closed monoidal $\V$-category $\A$,
generalizing the autonomous case;

\item we show the connection with the lax centre as well as the centre;
and,

\item we introduce the double $\D\A$ of a monoidal $\V$-category $\A$ and
some localizations of it, and relate these to Tambara modules.
\end{enumerate}

Our principal goal is to give conditions under which the centre and lax
centre of a $\V$-valued $\V$-functor monoidal $\V$-category is again such.
Some results in this direction can be found in~\cite{DaSt}.

For general monoidal $\A$, we construct a promonoidal $\V$-category
$\D\A$ with an equivalence $[\D\A,\V]\simeq \Tamb(\A)$. When $\A$ is
closed, we define when a Tambara module is (left) strong and show that
these constitute a $\V$-category ($\Tamb_{l s}(\A)$) $\Tamb_s(\A)$ which is
equivalent to the (lax) centre of $[\A,\V]$. We construct localizations
$\D_s\A$ and $\D_{l s}\A$ of $\D\A$ such that there are equivalences
$\Tamb_s(\A) \simeq [\D_s\A,\V]$ and $\Tamb_{l s}(\A) \simeq [\D_{l s}\A,\V]$.
When $\A$ is autonomous, every Tambara module is strong, which implies an
equivalence $\Z[\A,\V] \simeq [\D\A,\V]$. These results should be compared
with those of~\cite{DaSt} where the lax centre of $[\A,\V]$ is shown
generally to be a full sub-$\V$-category of a functor $\V$-category
$[\A_M,\V]$ which also becomes an equivalence $\Z[\A,\V] \simeq [\A_M,\V]$
when $\A$ is autonomous.

As we were completing this paper, Ignacio Lopez Franco sent us his
preprint~\cite{LF} which has some results in common with ours. As an example
for $\V$-modules of his general constructions on pseudomonoids, he is also
led to what we call the double monad.

\section{Centres and convolution}

We work with categories enriched in a base monoidal category $\V$ as used by
Kelly~\cite{KellyBook}. It is symmetric, closed, complete and cocomplete.

Let $\A$ denote a closed monoidal $\V$-category. We denote the tensor
product by $A \ox B$ and the unit by $I$ in the hope that this will cause no
confusion with the same symbols used for the base $\V$ itself. We have
$\V$-natural isomorphisms
\[
    \A(A,{}^B C) \cong \A(A \ox B,C) \cong \A(B,C^A)
\]
defined by evaluation and coevaluation morphisms
\begin{align*}
e_l &:{}^B C \ox B \ra C, & d_l &:A \ra {}^B (A \ox B), \\
e_r &:A \ox C^A \ra C, & d_r &:B \ra (A \ox B)^A.
\end{align*}
Consequently, there are canonical isomorphisms
\[
{}^{A \ox B} C \cong {}^A ({}^B C), \quad
C^{A \ox B} \cong (C^A)^B, \quad
({}^B C)^A \cong {}^B (C^A) \quad \mathrm{and} \quad
{}^I C \cong C \cong C^I
\]
which we write as if they were identifications just as we do with the
associativity and unit isomorphisms. We also write ${}^B C^A$ for
${}^B (C^A)$.

The Day convolution monoidal structure~\cite{Day} on the $\V$-category
$[\A,\V]$ of $\V$-functors from $\A$ to $\V$ consists of the tensor product
$F*G$ and unit $J$ defined by
\begin{align*}
(F*G)A &= \int^{U,V} \A(U \ox V,A) \ox FU \ox GV \\
       &\cong \int^V F ({}^V A) \ox GV \\
       &\cong \int^U FU \ox G(A^U)
\end{align*}
and
\[
    J A=\A(I,A).
\]
In particular,
\[
(F * \A(A,-)) B \cong F({}^A B) \qquad \mathrm{and} \qquad
(\A(A,-) * G) B \cong G(B^A).
\]

The centre of a monoidal category was defined in~\cite{JS91} and the lax
centre was defined, for example, in~\cite{DPS}. Since the representables are
dense in $[\A,\V]$, an object of the \emph{lax centre} $\Z_l[\A,\V]$ of
$[\A,\V]$ is a pair $(F,\theta)$ consisting of $F\in [\A,\V]$ and a
$\V$-natural family $\theta$ of morphisms
\[
    \theta_{A,B}:F({}^A B) \ra F(B^A)
\]
such that the diagrams
\[
\vcenter{\xygraph{{F({}^{A \ox B} C)}="s"
    :[r(3)] {F(C^{A \ox B})}="t" ^-{\theta_{A \ox B,C}}
 "s":[d] {F({}^A({}^B C))} _-=
    :[d(1.2)r(1.5)] {F({}^B C^A)} _-{\theta_{A,{}^B C}}
    :[u(1.2)r(1.5)] {F((C^A)^B)} _-{\theta_{B,C^A}}
    :"t" _-=}}
\qquad \text{and} \qquad
\vcenter{\xygraph{{F({}^I A)}="s"
    :[r(2)] {F(A^I)}="t" ^-{\theta_{I,A}}
 "s":[dr] {FA} _-=
    :"t" _-=}}
\]
commute. The hom object $\Z_l[\A,\V]((F,\theta),(G,\phi))$ is defined to be
the equalizer of two obvious morphisms out of $[\A,\V](F,G)$. The
\emph{centre} $\Z[\A,\V]$ of $[\A,\V]$ is the full sub-$\V$-category of
$\Z_l[\A,\V] $ consisting of those objects $(F,\theta)$ with $\theta$
invertible.

\section{Tambara modules}

Let $\A$ denote a monoidal $\V$-category. We do not need $\A$ to be closed
for the definition of Tambara module although we will require this
restriction again later.

A \emph{left Tambara module} on $\A$ is a $\V$-functor $T:\A^\op \ox \A
\ra \V$ together with a family of morphisms
\[
    \alpha_l(A,X,Y):T(X,Y) \ra T(A \ox X,A \ox Y)
\]
which are $\V$-natural in each of the objects $A$, $X$ and $Y$,
satisfying the two equations $\alpha_l(I,X,Y) = 1_{T(X,Y)}$ and
\[
\xygraph{{T(X,Y)}="s"
    :[d(1.4)r(2)] {T(A \ox A' \ox X,A \ox A' \ox Y).}="t"
         _-{\alpha_l(A \ox A',X,Y)}
 "s":[r(4)] {T(A' \ox X,A' \ox Y)} ^-{\alpha_l(A',X,Y)}
    :"t" ^-{\alpha_l(A,A' \ox X,A' \ox Y)}}
\]

Similarly, a \emph{right Tambara module} on $\A$ is a $\V$-functor
$T:\A^\op \ox \A \ra \V$ together with a family of morphisms
\[
\alpha_r(B,X,Y):T(X,Y) \ra T(X \ox B,Y \ox B)
\]
which are $\V$-natural in each of the objects $B$, $X$ and $Y$, satisfying
the two equations $\alpha_r(I,X,Y)=1_{T(X,Y)}$ and 
\[
\xygraph{{T(X,Y)}="s"
    :[d(1.4)r(2)] {T(X \ox B \ox B',Y \ox B \ox B').}="t"
         _-{\alpha_r(B \ox B',X,Y)}
 "s":[r(4)] {T(X \ox B,Y \ox B)} ^-{\alpha_r(B,X,Y)}
    :"t" ^-{\alpha_r(B',B \ox X,B \ox Y)}}
\]

A \emph{Tambara module}$ (T,\alpha)$ on $\A$ is a $\V$-functor $T:\A^\op
\ox \A \ra \V$ together with both left and right Tambara module structures
satisfying the ``bimodule'' compatibility condition
\[
\xygraph{{T(X,Y)}="s"
    :[r(5.5)] {T(A \ox X,A \ox Y)} ^-{\alpha_l(A,X,Y)}
    :[d(1.5)] {T(A \ox X \ox B,A\ox Y \ox B).}="t" ^-{\alpha_r(B,A\ox X,A\ox Y)}
 "s":[d(1.5)] {T(X \ox B,Y \ox B)} _-{\alpha_r(B,X,Y)}
    :"t" _-{\alpha_l(A,X \ox B,Y \ox B)}}
\]
The morphism defined to be the diagonal of the last square is denoted by
\[
    \alpha(A,B,X,Y):T(X,Y) \ra T(A \ox X \ox B, A \ox Y \ox B)
\]
and we can express a Tambara module structure purely in terms of this,
however, we need to refer to the left and right structures below.

\begin{proposition} \label{XRef-Proposition-812144421}
Suppose $\A$ is a monoidal $\V$-category and $T:\A^\op \ox \A \ra \V$ is a
$\V$-functor.
\begin{enumerate}[{\upshape (a)}]
\item If $\A$ is right closed, there is a bijection between $\V$-natural
families of morphisms
\[
    \alpha_l(A,X,Y):T(X,Y) \ra T(A \ox X,A \ox Y)
\]
and $\V$-natural families of morphisms
\[
    \beta_l(A,X,Y):T(X,Y^A) \ra T(A \ox X,Y).
\]

\item Under the bijection of (a), the family $\alpha_l$ is a left Tambara
structure if and only if the family $\beta_l$ satisfies the two equations 
$\beta_l(I,X,Y) = 1_{T(X,Y)}$ and 
\[
\xygraph{{T(X,Y^{A \ox A'})}="s"
    :[r(4.2)] {T(A \ox A' \ox X,Y)}="t" ^-{\beta_l(A \ox A',X,Y)}
 "s":[d(1.4)] {T(X,(Y^A)^{A'})} _-=
    :[r(4.2)] {T(A' \ox X,Y^A).} _-{\beta_l(A',X,Y^A)}
    :"t" _-{\beta_l(A,A'\ox X,Y)}}
\]

\item If $\A$ is left closed, there is a bijection between $\V$-natural
families of morphisms
\[
    \alpha_r(B,X,Y):T(X,Y) \ra T(X \ox B,Y \ox B)
\]
and $\V$-natural families of morphisms
\[
    \beta_r(B,X,Y):T(X,{}^B Y) \ra T(X \ox B,Y).
\]

\item Under the bijection of (c), the family $\alpha_r$ is a right Tambara
structure if and only if the family $\beta_r$ satisfies the two equations
$\beta_r(I,X,Y) = 1_T(X,Y)$ and
\[
\xygraph{{T(X,{}^{B \ox B'} Y)}="s"
    :[r(4.2)] {T(X \ox B \ox B',Y)}="t" ^-{\beta_r(B \ox B',X,Y)}
 "s":[d(1.4)] {T(X, {}^B ({}^{B'} Y))} _-=
    :[r(4.2)] {T(X \ox B,{}^{B'} Y).} _-{\beta_r(B,X,{}^{B'} Y)}
    :"t" _-{\beta_r(B',X \ox B,Y)}}
\]

\item If $\A$ is closed, the families $\alpha_l$ and $\alpha_r$ form a
Tambara module structure if and only if the families $\beta_l$ and
$\beta_r$, corresponding under (a) and (c), satisfy the condition
\[
\xygraph{{T(X,{}^B Y^A)}="s"
    :[r(4.4)] {T(A \ox X,{}^B Y)} ^-{\beta_l(A,X,{}^B Y)}
    :[d(1.4)] {T(A \ox X \ox B,Y).}="t" ^-{\beta_r(B,A \ox X,Y)}
 "s":[d(1.4)] {T(X \ox B,Y^A)} _-{\beta_r(B,X,Y^A)}
    :"t" _-{\beta_l(A,X \ox B,Y)}}
\]
\end{enumerate}
\end{proposition}

\begin{proof}
The bijection of (a) is defined by the formulas
\begin{multline*}
\beta_l(A,X,Y) = \Big(\xygraph{{T(X,Y^A)}
    :[r(3.4)] {T(A \ox X,A \ox Y^A)} ^-{\alpha_l(A,X,Y^A)}} \\
    \xygraph{:[r(2.4)] {T(A \ox X,Y)} ^-{T(A \ox X,e_r)}} \Big)
\end{multline*}
and
\begin{multline*}
\alpha_l(A,X,Y) = \Big(\xygraph{{T(X,Y)}
    :[r(3)] {T(X,(A \ox Y)^A)} ^-{T(X,d_r)}} \\
    \xygraph{:[r(3.0)] {T(A \ox X,A \ox Y)} ^-{\beta_l(A,X,A \ox Y)}}\Big).
\end{multline*}
That the processes are mutually inverse uses the adjunction identities on
the morphisms $e$ and $d$. The bijection of (c) is obtained dually by
reversing the tensor product. Translation of the conditions from the
$\alpha$ to the $\beta $ as required for (b), (d) and (e) is straightforward.
\end{proof}

A left (respectively, right) Tambara module $T$ on $\A$ will be called
\emph{strong} when the morphisms
\begin{align*}
    & \beta_l(A,X,Y):T(X,Y^A) \ra T(A \ox X,Y) \\
    & (\text{respectively,~}\ \beta_r(B,X,Y):T(X,{}^B Y) \ra T(X \ox B,Y))
\end{align*}
corresponding via Proposition~\ref{XRef-Proposition-812144421} to the left
(respectively, right) Tambara structure, are invertible. A Tambara module is
called \emph{left} (respectively, \emph{right}) \emph{strong} when it is
strong as a left (respectively, right) module and \emph{strong} when it is
both left and right strong. In particular, notice that the hom $\V$-functor
(= identity module) of $\A$ is a strong Tambara module.

\begin{proposition} \label{XRef-Proposition-9121222}
Suppose $\A$ is a monoidal $\V$-category and $T:\A^\op \ox \A \ra \V$ is a
$\V$-functor. If $\A$ is right (left) autonomous then every left (right)
Tambara module is strong.
\end{proposition}

\begin{proof}
If $A^*$ denotes a right dual for $A$ with unit $\eta:I \ra A^* \ox A$ then
an inverse for $\beta_l$ is defined by the composite
\[
\xygraph{{T(A \ox X,Y)}
    :[r(4.2)] {T(A^* \ox A \ox X,A^* \ox Y)} ^-{\alpha_l(A^*,A \ox X,Y)}
    :[r(3.5)] {T( X,A^* \ox Y)} ^-{T(\eta,1)}}.
\]
\end{proof}

Write $\LTamb(\A)$ for the $\V$-category whose objects are left Tambara
modules $T=(T,\alpha_l)$ and whose hom $\LTamb(\A) (T,T')$ in $\V$ is
defined to be the intersection over all $A$, $X$ and $Y$ of the equalizers
of the pairs of morphisms:
\[
\xygraph{{[\A^\op \ox \A,\V](T,T')}="s"
    :@<3pt>[r(5.8)] {\V(T(X,Y) ,T'(A \ox X,A \ox Y))}="t"
            ^-{\V(\alpha_l,1) \circ \pr_{A \ox X,A \ox Y}}
 "s":@<-3pt>"t" _-{\V(1,\alpha_l) \circ \pr_{X,Y}}}.
\]
Equivalently, we can define the hom as an intersection of equalizers
of pairs of morphisms:
\[
\xygraph{{[\A^\op \ox \A,\V](T,T')}="s"
    :@<3pt>[r(5.4)] {\V(T(X,Y^A) ,T'(A \ox X,Y))}="t"
            ^-{\V(\beta_l,1) \circ \pr_{A \ox X,Y}}
 "s":@<-3pt>"t" _-{\V(1,\beta_l) \circ \pr_{X,Y^A}}}.
\]
Composition is defined so that we have a $\V$-functor $\iota:\LTamb(\A) \ra
[\A^\op \ox \A,\V]$ which forgets the left module structure on $T$. In fact,
$\LTamb(\A)$ becomes a monoidal $\V$-category in such a way that the
forgetful $\V$-functor $\iota$ becomes strong monoidal. For this, the
monoidal structure on $[\A^\op \ox \A,\V]$ is the usual tensor product
(= composition) of endomodules:
\[
    (T \ox_\A T')(X,Y) = \int^Z T(X,Z) \ox T'(Z,Y).
\]
When $T$ and $T'$ are left Tambara modules, the left Tambara structure
\[
    (T \ox_\A T')(X,Y) \ra (T \ox_\A T')(A \ox X,A \ox Y)
\]
on $T \ox_\A T'$ is defined by taking its composite with the coprojection
$\copr_Z$ into the above coend to be the composite
\begin{multline*}
\xygraph{{T(X,Z) \ox T'(Z,Y)} 
    :[r(4.8)] {T(A \ox X,A \ox Z) \ox T'(A \ox Z,A \ox Y)}
            ^-{\alpha_l \ox \alpha_l}} \\
  \xygraph{:[r(3)] {(T\ox_\A T')(A \ox X,A \ox Y)} ^-{\copr_{A \ox Z}}}.
\end{multline*}
Similarly we obtain monoidal $\V$-categories $\RTamb(\A)$ and $\Tamb(\A)$ of
right Tambara and all Tambara modules on $\A$.

We write $\LTamb_s(\A)$ for the full sub-$\V$-category of $\LTamb(\A)$
consisting of the strong left Tambara modules. We write $\Tamb_{l s}(\A)$,
$\Tamb_{r s}(\A)$ and $\Tamb_s(\A)$ for the full sub-$\V$-categories of
$\Tamb(\A)$ consisting of the left strong, right strong and strong Tambara
modules respectively.

If $\A$ is autonomous then $\Tamb(\A) = \Tamb_{l s}(\A) = \Tamb_{r s}(\A) =
\Tamb_s(\A)$ by Proposition~\ref{XRef-Proposition-9121222}.

\section{The Cayley functor}

Consider a right closed monoidal $\V$-category $\A$. There is a
\emph{Cayley $\V$-functor}
\[
    \Upsilon:[\A,\V] \ra [\A^\op \ox \A,\V]
\]
defined as follows. To each object $F \in [\A,\V]$, define
$\Upsilon(F)=T_F$ by
\[
    T_F(X,Y) = F(Y^X).
\]
The effect $\Upsilon_{F,G}:[\A,\V](F,G) \ra [\A^\op \ox \A,\V](T_F,T_G)$
of $\Upsilon$ on homs is defined by taking its composite with the
projection
\[
    \pr_{X,Y}:[\A^\op \ox \A,\V](T_F,T_G) \ra \V( F(Y^X) ,G(Y^X))
\]
to be the projection
\[
    \pr_{Y^X}:[\A,\V](F,G) \ra \V(F(Y^X),G(Y^X)).
\]

\begin{proposition}
The Cayley $\V$-functor $\Upsilon$ is strong monoidal; it takes Day
convolution to composition of endomodules.
\end{proposition}

\begin{proof}
We have the calculation:
\begin{align*}
(\Upsilon(F) \ox_\A \Upsilon(G))(X,Y)
    &= \int^Z \Upsilon(F)(X,Z) \ox \Upsilon(G)(Z,Y) \\
    &= \int^Z F(Z^X) \ox G(Y^Z) \\
    &\cong \int^{Z,U,V} \A(U,Z^X) \ox FU \ox \A(V,Y^Z) \ox GV \\
    &\cong \int^{Z,U,V} \A(X \ox U,Z) \ox FU \ox \A(Z \ox V,Y) \ox GV \\
    &\cong \int^{U,V} \A(X \ox U \ox V,Y) \ox FU \ox GV \\
    &\cong \int^{U,V} \A(U \ox V,Y^X) \ox FU \ox GV \\
    &\cong \Upsilon(F*G)(X,Y),
\end{align*}
and of course $\Upsilon(\A(I,-))(X,Y) = \A(I,Y^X) \cong \A(X,Y)$.
\end{proof}

In fact, $\Upsilon$ lands in the left Tambara modules by defining, for each
$F\in [\A,\V]$, the structure
\[
\alpha_l(A,X,Y) =\Big(
    \xymatrix@C=10ex{F(Y^X) \ar[r]^-{F((d_r)^X)} & F((A \ox Y)^{A \ox X})}
    \Big)
\]
on $T_F$. It is helpful to observe that the $\beta_l$ corresponding to this
$\alpha_l$ (via Proposition~\ref{XRef-Proposition-812144421}) is given by
the identity
\[
\beta_l(A,X,Y) = \Big(
    \xymatrix@C=5ex{F(Y^{A \ox X}) \ar[r]^-1 & F(Y^{A \ox X})}\Big),
\]
showing that $T_F$ becomes a strong left module. To see that there is a
$\V$-functor $\hat{\Upsilon}:[\A,\V] \ra \LTamb_s(\A)$ satisfying $\iota
\circ \hat{\Upsilon}=\Upsilon $, we merely observe that
\[
\pr_{A \ox X,Y} \circ \Upsilon_{F,G} = \pr_{Y^{A \ox X}} = \pr_{(Y^A)^X}
=\pr_{X,Y^A} \circ \Upsilon_{F,G}.
\]

\begin{proposition} \label{XRef-Proposition-83065512}
If $\A$ is a right closed monoidal $\V$-category then the $\V$-functor
$\hat{\Upsilon}:[\A,\V] \ra \LTamb_s(\A) $ is an equivalence.
\end{proposition}

\begin{proof}
Define $\zeta:\LTamb(\A)(T_F,T_G) \ra [\A,\V](F,G)$ by $\pr_Y \circ \zeta =
\pr_{I,Y} \circ \iota_{T_F,T_G}$. Then
\[
\pr_Y \circ \zeta \circ \hat{\Upsilon}_{F,G}
    = \pr_{I,Y} \circ \iota_{T_F,T_G} \circ \hat{\Upsilon}_{F,G}
    = \pr_{I,Y} \circ \Upsilon_{F,G}
    = \pr_Y
\]
and
\begin{align*}
\pr_{X,Y} \circ \iota _{T_F,T_G} \circ \hat{\Upsilon}_{F,G} \circ \zeta
    &= \pr_{X,Y} \circ \Upsilon_{F,G} \circ \zeta \\
    &= \pr_{Y^X} \circ \zeta \\
    &= \pr_{I,Y^X} \circ \iota_{T_F,T_G} \\
    &= \pr_{X,Y} \circ \iota_{T_F,T_G}.
\end{align*}
It follows that $\zeta $ is the inverse of $\hat{\Upsilon}_{F,G}$, so that
$\hat{\Upsilon}$ is fully faithful. To see that $\hat{\Upsilon}$ is
essentially surjective on objects, take a strong left module $T$. Put
$FY=T(I,Y) $ as a $\V$-functor in $Y$. Then the isomorphism $\beta_l(X,I,Y)$
yields
\[
    T_F(X,Y) = F(Y^X) = T(I,Y^X) \cong T(X,Y);
\]
so $\hat{\Upsilon}(F) \cong T$.
\end{proof}

Now suppose we have an object $(F,\theta)$ of the lax centre $\Z_l[\A,\V]$
of $[\A,\V]$. Then $T_F$ becomes a right Tambara module by defining
\[
\alpha_r(B,X,Y) = \Big(\xygraph{{F(Y^X)}
    :[r(2.8)] {F({}^{B}(Y \ox B)^X)} ^-{F((d_l)^X)}
    :[r(3.2)] {F(Y \ox B)^{X \ox B}} ^-{\theta_{B,(Y \ox B)^X}}}\Big).
\]
If $\A$ is left closed, the $\beta_r$ corresponding to this $\alpha_r$
(via Proposition~\ref{XRef-Proposition-812144421}) is defined by
\[
    \beta_r(B,X,Y) =\Big(\xymatrix@C=7ex{
    F({}^B Y^X) \ar[r]^-{\theta_{B,Y^X}} & F(Y^{X \ox B})}\Big).
\]
It is easy to see that, in this way, $T_{F}=\hat{\Upsilon}(F)$ actually
becomes a (two-sided) Tambara module which we write as
$\hat{\Upsilon}(F,\theta)$, and we have a $\V$-functor
\[
    \hat{\Upsilon}:\Z_l[\A,\V] \ra \Tamb_{l s}(\A).
\]

\begin{proposition}
If $\A$ is a closed monoidal $\V$-category then the $\V$-functor
\[
    \hat{\Upsilon}:\Z_l[\A,\V] \ra \Tamb_{l s}(\A)
\]
is an equivalence which restricts to an equivalence
\[
    \hat{\Upsilon}:\Z[\A,\V] \ra \Tamb_s(\A).
\]
\end{proposition}

\begin{proof}
The proof of full faithfulness proceeds along the lines of the beginning of
the proof of Proposition~\ref{XRef-Proposition-83065512}. For essential
surjectivity on objects, take a left strong Tambara module $(T,\alpha)$.
Then $\beta_l(A,X,Y):T(X,Y^A) \ra T(A \ox X,Y)$ is invertible. Define the
$\V$-functor $F:\A \ra \V$ by $FX = T(I,X)$ as in the proof of
Proposition~\ref{XRef-Proposition-83065512}, and define
$\theta_{A,Y}:F({}^A Y) \ra F(Y^A)$ to be the composite
\[
\xygraph{{T(I,{}^A Y)}
    :[r(2.8)] {T(A,Y)} ^-{\beta_r(A,I,Y)}
    :[r(3)] {T(I,{}Y^A)} ^-{\beta_l(A,I,Y)^{-1}}}.
\]
This is easily seen to yield an object $(F,\theta)$ of the lax centre
$\Z_l[\A,\V] $ with $\hat{\Upsilon}(F,\theta) \cong T_F$. Thus we have the
first equivalence. Clearly $\theta$ is invertible if and only if $\beta_r$
is; the second equivalence follows.
\end{proof}

\section{The double monad}

Tambara modules are actually Eilenberg-Moore coalgebras for a fairly
obvious comonad on $[\A^\op \ox \A,\V]$. We begin by looking at the case
of left modules.

Let $\Theta _{l}:[\A^\op \ox \A,\V] \ra [\A^\op \ox \A,\V]$ be the
$\V$-functor defined by the end
\[
    \Theta_l(T)(X,Y) = \int_A T(A \ox X,A \ox Y).
\]
There is a $\V$-natural family $\epsilon_{T}:\Theta_l(T) \ra T$ defined
by the projections
\[
    \pr_I:\int_A T(A \ox X,A \ox Y) \ra T(X,Y).
\]
There is a $\V$-natural family $\delta_T:\Theta_l(T) \ra \Theta_l(\Theta_l(T))$
defined by taking its composite with the projection
\[
    \pr_{B,C}:\int_{B,C} T(B \ox C \ox X,B \ox C \ox Y) \ra
    T(B \ox C \ox X,B \ox C \ox Y)
\]
to be the projection
\[
  \pr_{B \ox C}:\int_A T(A \ox X, A \ox Y) \ra T(B \ox C \ox X, B \ox C \ox Y).
\]
It is now easily checked that $\Theta_l= (\Theta_l,\delta,\epsilon)$ is a
comonad on $[\A^\op \ox \A,\V]$.

There are also a comonad $\Theta_r$ on $[\A^\op \ox \A,\V]$, a distributive
law $\Theta_r \Theta_l \cong \Theta_l \Theta_r$, and a comonad $\Theta =
\Theta_r \Theta_l$:
\[
    \Theta_r(T)(X,Y) = \int_B T(X \ox B, Y \ox B)
\]
and
\[
    \Theta(T)(X,Y) = \int_{A,B} T(A \ox X \ox B, A \ox Y \ox B).
\]
We can easily identify the $\V$-categories of Eilenberg-Moore coalgebras
for these three comonads.

\begin{proposition}\label{XRef-Proposition-96133952}
There are isomorphisms of $\V$-categories
\begin{itemize}
\item $[\A^\op \ox \A,\V]^{\Theta_l} \cong \LTamb(\A)$,
\item $[\A^\op \ox \A,\V]^{\Theta_r} \cong \RTamb(\A)$, and
\item $[\A^\op \ox \A,\V]^\Theta \cong \Tamb(\A)$.
\end{itemize}
\end{proposition}

In fact, $\Theta_l$, $\Theta_r$ and $\Theta$ are all monoidal comonads on
$[\A^\op \ox \A,\V]$. For example, the structure on $\Theta_l$ is provided
by the $\V$-natural transformations $\Theta_l(T) \ox_\A \Theta_l(T')$ $\ra$
$\Theta_l(T \ox_\A T')$ and $\A(-,-) \ra \Theta_l(\A(-,-))$ with components
\begin{equation}\label{XRef-Equation-96131613}
    \int^Z\!\!\! \int_A T(A \ox X, A \ox Z) \ox\!\! \int _{B}T'(B \ox X,B\ox Z)
    \ra \int_C \int^U\!\!\! T(C \ox X,U) \ox T'(U,C \ox Y)
\end{equation}
and
\begin{equation} \label{XRef-Equation-96132817}
    \A(X,Y) \ra \int_A \A(A \ox X, A \ox Y)
\end{equation}
defined as follows. The morphism~\eqref{XRef-Equation-96131613} is
determined by its precomposite with the coprojection $\copr_Z$ and
postcomposite with the projection $\pr_C$; the result is defined to be the
composite
\begin{align*}
\int_A T(A \ox X, A \ox Z) \ox & \int_B T'(B\ox X,B\ox Z) \\
    & \xymatrix@C=12ex{\ar[r]^-{\pr_C \>\ox\> \pr_C}&}
            {T(C \ox X,C \ox Z) \ox T'(C \ox Z,C \ox Y)} \\
    & \xymatrix@C=12ex{\ar[r]^-{\copr_{C \ox Z}}&}
            \int^U T(C \ox X,U) \ox T'(U,C \ox Y)~.
\end{align*}
The morphism~\eqref{XRef-Equation-96132817} is simply the coprojection
$\copr_I$. It follows that $[\A^\op \ox \A,\V]^{\Theta_l}$ becomes monoidal
with the underlying functor becoming strong monoidal; see~\cite{Moerd02}
and~\cite{McC02}. Clearly we have:

\begin{proposition}
The isomorphisms of Proposition~\ref{XRef-Proposition-96133952} are monoidal.
\end{proposition}

The next thing to observe is that $\Theta_l$, $\Theta_r$ and $\Theta$ all
have left adjoints $\Phi_l$, $\Phi_r$ and $\Phi$ which therefore become
opmonoidal monads whose $\V$-categories of Eilenberg-Moore algebras are
monoidally isomorphic to $\LTamb(\A)$, $\RTamb(\A)$ and $\Tamb(\A)$,
respectively. Straightforward applications of the Yoneda Lemma, show that
the formulas for these adjoints are
\begin{align*}
& \Phi_l(S)(U,V) = \int^{A,X,Y} \A(U,A \ox X) \ox \A(A \ox Y,V) \ox S(X,Y), \\
& \Phi_r(S)(U,V) = \int^{B,X,Y} \A(U,X \ox B) \ox \A(Y \ox B,V) \ox S(X,Y),
    \quad \text{and} \\
& \Phi(S)(U,V) = \int^{A,B,X,Y} \A(U,A \ox X \ox B) \ox \A(A \ox Y \ox B,V)
    \ox S(X,Y).
\end{align*}

Recall that left adjoint $\V$-functors $\Psi:[\X^\op,\V] \ra [\Y^\op,\V]$
are equivalent to $\V$-functors $\check{\Psi}: \Y^\op \ox \X \ra \V$, which
are also called modules $\check{\Psi}:\X \mra \Y$ from $\X$ to $\Y$. The
equivalence is defined by:
\[
    \check{\Psi}(Y,X) = \Psi(\X(-,X))(Y)
\]
and
\[
    \Psi(M)(Y) = \int^X \check{\Psi}(Y,X) \ox M(X).
\]

It follows that $\Phi_l$, $\Phi_r$ and $\Phi$ determine monads
$\check{\Phi}_l$, $\check{\Phi}_r$ and $\check{\Phi}$ on $\A^\op \ox \A$
in the bicategory $\VMod$. The formulas are:
\begin{align*}
& \check{\Phi}_l (X,Y,U,V) = \int^A \A(U,A \ox X) \ox \A(A \ox Y,V), \\
& \check{\Phi}_r (X,Y,U,V) = \int^B \A(U,X \ox B) \ox \A(Y \ox B,V),
    \quad \text{and} \\
& \check{\Phi}(X,Y,U,V) = \int^{A,B} \A(U,A \ox X \ox B) \ox \A(A \ox Y\ox B,V).
\end{align*}

\section{Doubles}\label{XRef-Section-831214024}

The bicategory $\VMod$ admits the Kleisli construction for monads. Write
$\D_l\A$, $\D_r\A$ and $\D\A$ for the Kleisli $\V$-categories for the
monads $\check{\Phi}_l$, $\check{\Phi}_r$ and $\check{\Phi}$ on
$\A^\op \ox \A$ in the bicategory $\VMod$. We call them the \emph{left
double}, \emph{right double} and \emph{double} of the monoidal
$\V$-category $\A$. They all have the same objects as $\A^\op \ox \A$. The
homs are defined by
\begin{align*}
& \D_l\A((X,Y), (U,V)) = \int^A \A(U,A \ox X) \ox \A(A \ox Y,V), \\
& \D_r\A((X,Y), (U,V)) = \int^B \A(U,X \ox B) \ox \A(Y \ox B,V),
    \quad \text{and} \\
& \D\A((X,Y), (U,V)) = \int^{A,B} \A(U,A \ox X \ox B) \ox \A(A \ox Y \ox B,V).
\end{align*}

\begin{proposition} \label{XRef-Proposition-98174931}
There are canonical equivalences of $\V$-categories:
\begin{itemize}
\item $\Xi_l :\LTamb(\A) \simeq [\D_l\A,\V]$,
\item $\Xi_r :\RTamb(\A) \simeq [\D_r\A,\V]$, and
\item $\Xi   :\Tamb(\A) \simeq [\D\A,\V]$.
\end{itemize}
\end{proposition}

It follows from the main result of Day~\cite{Day} that these doubles
$\D_l\A$, $\D_r\A$ and $\D\A$ all admit promonoidal structures $(P_l,J_l)$,
$(P_r,J_r)$ and $(P,J)$ for which the equivalences in
Proposition~\ref{XRef-Proposition-98174931} become monoidal when the
right-hand sides are given the corresponding convolution structures. For
example, we calculate that $P_l$ and $J_l$ are as follows:
\begin{align*}
P_l((X,Y),& (U,V);(H,K)) = (\D_l\A((X,Y),-) \ox_\A \D_l\A((U,V),-))(H,K) \\
    &= \int^{Z,A,B} \A(H,A \ox X) \ox \A(A \ox Y,Z) \ox \A(Z,B \ox U) \ox
        \A(B \ox V,K) \\
    &= \int^{A,B} \A(H,A \ox X) \ox \A(A \ox Y,B \ox U) \ox \A(B \ox V,K)
\end{align*}
and
\[
    J_l(H,K) = \A(H,K).
\]

Furthermore, there are some special morphisms that exist in these doubles
$\D_l\A$, $\D_r\A$ and $\D\A$. Let $\tilde{\alpha}_l:(X,Y) \ra
(A \ox X,A \ox Y)$ denote the morphism in $\D_l\A$ defined by the composite
\begin{multline*}
\xygraph{{I} :[r(4.2)] {\A(A \ox X,A \ox X) \ox \A(A \ox Y,A \ox Y)}
                   ^-{j_{A \ox X} \ox j_{A \ox Y}}} \\
    \xygraph{:[r(3.1)] {\D_l\A((X,Y),(A \ox X, A \ox Y))} ^-{\copr_A}}.
\end{multline*}
The $\V$-functor $\Xi_l$ has the property that $\Xi_l(T,\alpha_l)(X,Y) =
T(X,Y)$ and $\Xi_l(T,\alpha_l)(\tilde{\alpha}_l) = \alpha_l$. When $\A$ is
right closed, we let $\tilde{\beta}_l:(X,Y^A) \ra (A \ox X,Y)$ denote the
morphism in $\D_l\A$ defined by the composite
\begin{multline*}
\xygraph{{I}
:[r(3.8)] {\A(A \ox X,A \ox X) \ox \A(A \ox Y^A,Y)} ^-{j_{A \ox X} \ox e_r}} \\
  \xygraph{:[r(3)] {\D_l\A((X,Y^A),(A \ox X,Y))} ^-{\copr_A}}.
\end{multline*}
Then $\Xi_l(T,\alpha_l)(\tilde{\beta}_l) = \beta_l$.

Similarly, we have the morphism $\tilde{\alpha}_r:(X,Y) \ra (X \ox B,Y \ox B)$
in $\D_r\A$, and also, when $\A$ is left closed, the morphism
$\tilde{\beta}_r:(X,{}^B Y) \ra (X \ox B,Y)$.

There are $\V$-functors $\D_l\A \ra \D\A \la \D_r\A$ which are the identity
functions on objects and are defined on homs using projections with $B=I$
for the left leg and the projections $A=I$ for the second leg. In this way,
the morphisms $\tilde{\alpha}_l$ and $\tilde{\alpha}_r$ can be regarded also
as morphisms of $\D\A$. Under closedness assumptions, the morphisms
$\tilde{\beta}_l$ and $\tilde{\beta}_r$ can also be regarded as morphisms of
$\D\A$.

Let $\Sigma_l$ denote the set of morphisms $\tilde{\beta}_l:(X,Y^A) \ra
(A \ox X,Y)$, let $\Sigma_r$ denote the set of morphisms $\tilde{\beta}_r:
(X,{}^B Y) \ra (X \ox B,Y)$, and let $\Sigma$ denote the set of morphisms
$\Sigma = \Sigma_l \cup \Sigma_r$. Under appropriate closedness assumptions
on $\A$, we can form various $\V$-categories of fractions such as:
\begin{itemize}
\item $\mathrm{L}\D\A = \D_l\A[\Sigma_l^{-1}]$~~and~
    $\mathrm{R}\D\A  =\D_r\A[\Sigma_r^{-1}]$,
\item $\D_{l s}\A = \D\A[\Sigma_l^{-1}]$~~and~
    $\D_{r s}\A = \D\A[\Sigma_r^{-1}]$,~~and
\item $\D_s\A = \D\A[\Sigma^{-1}]$.
\end{itemize}
The following result is now automatic.

\begin{theorem}\label{XRef-Theorem-831214130}
For a closed monoidal $\V$-category $\A$, there are equivalences of
$\V$-categories:
\begin{itemize}
\item $[\mathrm{L}\D\A,\V] \simeq \LTamb_s(\A) \simeq [\A,\V]$,
\item $[\D_{l s}\A,\V] \simeq \Tamb_{l s}(\A) \simeq \Z_l[\A,\V]$, and
\item $[\D_s\A,\V] \simeq \Tamb_s(\A) \simeq \Z[\A,\V]$.
\end{itemize}
\end{theorem}

The first equivalence of Theorem~\ref{XRef-Theorem-831214130} implies
that $\mathrm{L}\D\A$ and $\A$ are Morita equivalent. This begs the question
of whether there is a $\V$-functor relating them more directly. Indeed there
is. We have a $\V$-functor
\[
    \Pi:\D_l\A \ra \A
\]
defined on objects by $\Pi(X,Y)=Y^X$ and by defining the effect
\[
    \Pi:\D_l\A((X,Y),(U,V)) \ra \A(Y^X,V^U)
\]
on hom objects to have its composite with the $A$-coprojection equal
to the composite
\begin{align*}
 \A(U,A \ox X) & \ox \A(A \ox Y,V) \\
& \xymatrix@C=14ex{\ar[r]^-{V^{(-)} \ox (-)^{A \ox X}} &
    \A(V^{A \ox X}, V^U) \ox \A((A \ox Y)^{A \ox X}, V^{A \ox X})} \\
& \xymatrix@C=14ex{\ar[r]^-{\textrm{composition}} &
    \A((A \ox Y)^{A \ox X},V^U)} \\
& \xymatrix@C=14ex{\ar[r]^-{\A((d_r)^{X},V^U)} & \A(Y^X,V^U)}.
\end{align*}
It is easy to see that $\Pi$ takes the morphisms $\tilde{\beta}_l:(X,Y^A)
\ra (A \ox X,Y)$ to isomorphisms. So $\Pi$ induces a $\V$-functor
\[
    \hat{\Pi}:\mathrm{L}\D_l\A \ra \A;
\]
this induces the first equivalence of Theorem~\ref{XRef-Theorem-831214130}.

For closed monoidal $\A$, the second and third equivalences of
Theorem~\ref{XRef-Theorem-831214130} show that both the lax centre and the
centre of the convolution monoidal $\V$-category $[\A,\V]$ are again functor
$\V$-categories $[\D_{l s}\A,\V]$ and $[\D_s\A,\V]$. Since $\Z_l[\A,\V]$
and $\Z[\A,\V]$ are monoidal with the tensor products colimit preserving in
each variable, using the correspondence in~\cite{Day}, there are lax braided
and braided promonoidal structures on $\D_{l s}\A$ and $\D_s\A$ which are
such that $[\D_{l s}\A,\V]$ and $[\D_s\A,\V]$ become closed monoidal under
convolution, and the equivalences of Theorem~\ref{XRef-Theorem-831214130}
become lax braided and braided monoidal equivalences.

\begin{remark} \label{XRef-Remark-91172928}
We are grateful to Brian Day for pointing out that the $\V$-category
$\A_M$ appearing in~\cite{DaSt} is equivalent to the full sub-$\V$-category
of $\D\A$ consisting of the objects of the form $(I,Y)$.

He also pointed out that a consequence of
Theorem~\ref{XRef-Theorem-831214130} is that the centre of $\V$ as a
$\V$-category is equivalent to $\V$ itself. This also can be seen directly
by using the $\V$-naturality in $X$ of the centre structure $u_X:A \ox X
\ra X \ox A$ on an object $A$ of $\V$, and the fact that $u_I=1$, to deduce
that $u_X = c_{A,X}$ (the symmetry of $\V$). Generally, the centre of $\V$
as a monoidal $\Set$-category is not equivalent to $\V$.
\end{remark}


\end{document}